\documentclass[11pt]{article}

\usepackage[a4paper, total={6.5in, 9in}]{geometry}

\def\bfs{{\mathbf s}}
\usepackage{amsmath,mathrsfs}
\def\MR{MR}

\usepackage{shapepar,bm}
\usepackage[dvips,arrow,matrix,ps,color,line]{xy}
\usepackage{theorem,amsfonts,amssymb,amscd}
\usepackage{geometry}
\usepackage{makeidx}
\usepackage{stmaryrd}

\providecommand{\bysame}{\leavevmode\hbox to3em{\hrulefill}\thinspace}
\providecommand{\MR}{\relax\ifhmode\unskip\space\fi MR }

\providecommand{\href}[2]{#2}

\usepackage{graphicx}
\usepackage{graphics}

\usepackage{epstopdf}

\usepackage{amssymb,graphics,hyperref}

\hypersetup{
    colorlinks = true,
    linkcolor = blue,
    anchorcolor = red,
   citecolor = blue,
    filecolor = red,
    pagecolor = red,
    urlcolor = blue}

\newcommand{\mathsym}[1]{{}}

\usepackage[usenames]{color}
\definecolor{MyLightMagenta}{cmyk}{0.1,0.8,0,0.1}
\definecolor{MyDarkBlue}{rgb}{0.1,0,0.3}

\makeindex

\def\wb{[\bfb]}
\def\wbet{[\bm\beta]}

\def\NN{\mathbb N}
\def\DJKM{{Date, Jimbo, Kashiwara and Miwa}}

\def\bfz{{\mathbf z}}

\def\ovsig{\overline{\sigma}}

\def\bbeta{{\bm\beta}}
\newcommand{\llb}{\llbracket}
\newcommand{\rrb}{\rrbracket}
\def\bfb{{\mathbf b}}

\def\ssetminus{\hskip-3pt\setminus\hskip-2pt}

\def\bfu{{\mathbf u}}
\def\bfv{{\mathbf v}}
\def\bfw{{\mathbf w}}

\def\d{\partial}

\def\Vcal{\mathcal V}

\def\bfx{{\mathbf x}}

\def\ZZ{\mathbb Z}

\def\Ecal{{\mathcal E}}

\def\QQ{\mathbb Q}

\def\cocoa{{\hbox{\rm C\kern-.13em o\kern-.07em C\kern-.13em o\kern-.15em A}}}

\def\Dcal{\mathcal D}

\def\Fcal{{\mathcal F}}

\def\hom{{\mathrm{Hom}}}

\def\bfu{{\bf u}}

\def\End{\mathrm{End}}

\def\blamb{{\bm \lambda}}
\def\bmu{{\bm \mu}}
\def\bnu{{\bm \nu}}
\def\Pcal{{\mathcal P}}

\def\w2M{\bigwedge^2M}

\def\w{\wedge }
\def\bw{\bigwedge }

\protect
\protect
\protect
\protect

\protect
\protect

\def\sra{\rightarrow}

\def\proof{\noindent{\bf Proof.}\,\,}
\def\qed{{\hfill\vrule height4pt width4pt depth0pt}\medskip}
\def\be{\begin{equation}}
\def\ee{\end{equation}}
\def\bclm{\begin{claim}}
\def\eclm{\end{claim}}
\def\beqn{\begin{eqnarray}}
\def\eeqn{\end{eqnarray}}
\def\beqn*{\begin{eqnarray*}}
\def\eeqn*{\end{eqnarray*}}
\def\ovPc{\overline{\Pcal}}

\theoremstyle{change}
\theorembodyfont{\rmfamily}

\newtheorem{claim}{}[section]

\global
\global

\def\no@breaks#1{{\def\\{ \ignorespaces}#1}}    

  \makeatletter
\def\cleardoublepage{\clearpage\if@twoside \ifodd\c@page\else
\hbox{} \thispagestyle{empty}
\newpage
\if@twocolumn\hbox{}\newpage\fi\fi\fi} \makeatother

\usepackage{eso-pic,graphicx}
\makeatletter
\newcommand\BackgroundPicture[2]{%
  \setlength{\unitlength}{1pt}%
  default \put(0,\strip@pt\paperheight){%
  \parbox[t][\paperheight]{\paperwidth}{%
    \vfill
     \centering \includegraphics[angle=#2, width=15cm, height=15cm,  bb=0 0 150 150]{#1}
    \vfill
}}} %
\makeatother

\date{}

\title{Bosonic and Fermionic  Representations of Endomorphisms of Exterior Algebras\\
\author{{\sc Ommolbanin Behzad, Letterio Gatto}
\thanks{Work sponsored by  Finanziamento 
Diffuso della Ricerca, no. 53$\_$RBA17GATLET del Politecnico di Torino; 
\, Progetto di Eccellenza\, Dipartimento\, di\, Scienze\, Matematiche, 2018--2022 no.
E11G18000350001, INDAM-GNSAGA e PRIN "Geometria delle Variet\`a Algebriche''.
\smallskip 
\newline  ${}$
\,\,\,\,\,\,\,{\em Keywords and Phrases:} Schubert Derivations on the fermionic Fock space, Vertex Operators  on Exterior Algebras,  Bosonic and Fermionic Representations by Date-Jimbo-Kashiwara-Miwa, Symmetric Functions.
\smallskip
\newline ${}$ \,\,\,\,\,\, {\bf 2020 MSC:} 14M15, 15A75, 05E05, 17B69. \,}}}

\begin{document}

\maketitle

\begin{abstract} \noindent We describe the fermionic and bosonic Fock  representation  of the Lie super-algebra of   endomorphisms
of the exterior algebra of the $\QQ$-vector space of infinite countable dimension, vanishing at all but finitely many basis elements. We achieve the goal by exploiting the extension of the Schubert derivations on the fermionic Fock space.
\end{abstract}

\section*{Introduction} 

\claim{\bf The Goal.} Let  $B:=
\QQ[\bfx]$ be the polynomial ring in  the infinitely many indeterminates 
$\bfx:=(x_1,x_2,\ldots)$. The purpose of this paper is to further enhance 
a classical but fundamental result by \DJKM, which describes  the 
polynomial algebra  
$B(\xi):=B \otimes_\QQ\QQ[\xi^{-1},\xi]$ (the {\em bosonic Fock space}) 
as a 
representation of the Lie algebra 
\begin{center}
$gl_\infty(\QQ):=\{(a_{ij})_{i,j\in\ZZ}\,|\,$ all but a finite number of the $a_{ij}\in\QQ$ are zero$\}$. 
\end{center}
For the convenience of our exposition, we slightly change the notation for $gl_\infty(\QQ)$ as follows. Let us consider the $\QQ$-vector space $\Vcal:=\bigoplus_{i\in
\ZZ}\QQ\cdot b_i$, whose basis $\bfb=(b_i)_{i\in\ZZ}$ is parametrised by the integers,  along with its restricted dual $\Vcal^*:=\bigoplus_{i\in\ZZ}\QQ
\cdot \beta_j$,  with basis $\bbeta=(\beta_i)_{i\in\ZZ}$, where $\beta_j \in \hom_\QQ(\Vcal, \QQ)$  is the unique 
linear form such that $\beta_j (b_i)=\delta_{ji}$. 
 In this way, the algebra $gl_\infty(\QQ)$ gets
identified with 
\be
gl(\Vcal)=\Vcal\otimes \Vcal^*=\bigoplus_{i,j\in\ZZ}\QQ\cdot b_i\otimes
\beta_j.\label{eq0:basbvr}
\ee
In the contribution \cite{DJKM01} (see also \cite{jimbomiwa}) \DJKM\, 
compute 
the action on $B(\xi)$ of the generating function 
\be
\Ecal(z,w)=\sum_{i,j}
b_i\otimes\beta_jz^iw^{-j}\label{eq0:gfdjkm}
\ee 
of the basis as in \eqref{eq0:basbvr}, and obtain
their celebrated {\em bosonic vertex operator representation} of $gl(\Vcal)$.
Our main result is  the ultimate extension of the DJKM formula, concerned with  the
fermionic and the bosonic vertex representation of the Lie super--algebra 
$gl(\bw\Vcal):=\bw\Vcal\otimes\bw\Vcal^*$, where  $\bw\Vcal$ and 
$\bw\Vcal^*$ denote, as usual, the exterior algebra of $\Vcal$ and $\Vcal^*$ 
respectively. Our formula includes, and admits as a particular 
case, the DJKM one, as displayed, e.g., in \cite[Proposition 5.2]{KacRaRoz} or \cite[Section 1]{jimbomiwa}.

\claim{} Let $\Pcal$ be the set of all partitions (non 
increasing sequences of non negative integers all zero but finitely 
many). The {\em fermionic Fock space} is a $\ZZ$-graded 
vector space $\Fcal:=\bigoplus_{m\in\ZZ}\Fcal_m$ which, like $B(\xi)$, the bosonic 
one, possesses a basis $\wb_{m+\blamb}$ parametrised by $\ZZ
\times \Pcal$. More than that, it is essential, from our point of view, to think of $\Fcal$ as a $B(\xi)$-module of rank $1$ 
generated 
by 
$
\wb_0:=b_0\w b_{-1}\w b_{-2}\w\cdots,$
 such that 
\be
\xi^mS_\blamb(\bfx)\wb_0=\wb_{m+\blamb}=b_{m+\lambda_1}\w\cdots\w b_{m-r
+1+\lambda_r}\w b_{m-r}\w b_{m-r-1}\w\cdots,
\label{eq0:pbfc}
\ee
where  $S_\blamb(\bfx)$ denotes the Schur polynomial 
associated to the partition $\blamb$ and to the sequence $\bfx$.
Equality~\eqref{eq0:pbfc} can be  understood either as a  Giambelli's\,\, 
formula for Schubert Calculus on infinite Grassmannian (see \cite{gln}) 
or like a Jacobi-Trudy formula. To follow more closely the reference 
\cite[Theorem 6.1]{KacRaRoz}, and being more adherent to the subject of the paper, we 
call  \eqref{eq0:pbfc}  the
{\em Boson-Fermion correspondence}. Our starting point is the obvious remark that $\bw
\Vcal$ is a (irreducible)  representation of the Lie super--algebra $gl(\bw
\Vcal)$ of all endomorphisms vanishing at 
all  basis elements but finitely many of the exterior algebra. An explicit 
generating function encoding  the $gl(\bw\Vcal)$-module structure of $\bw\Vcal$ has 
already been proposed in \cite{BeCoGaVi}, where  the  vertex 
operators shaping the boson-fermion correspondence  spontaneously arise in all 
their splendor, regardless of the  more classical framework. In addition, as
noticed in \cite{SDIWP}, little effort is needed  to 
extend the $\bw\Vcal$-representation to $\Fcal$, mainly because the latter is a module over the former. This reflects in the fact that each degree 
$\Fcal_m$ of $\Fcal$, as suggested in formula \eqref{eq0:pbfc},  can be thought 
of as a semi-infinite 
exterior power. Finally, one  just pulls back on $B(\xi)$ the $\Fcal$ representation of $gl(\bw\Vcal)$, invoking the boson-fermion correspondence. The program demands, however,  to identify a basis of $\bw\Vcal\otimes \bw\Vcal^*$ suited to get a convenient
generalisation of the DJKM generating function \eqref{eq0:djkmmm}. Last, but not the least,  one  is left to determine explicitly its action on $\bw\Vcal$. This  is the  point that, as in our previous contribution,  the flexible formalism of Schubert derivations (a distinguished kind of Hasse-Schmidt derivation on an exterior algebra), extended to $\Fcal$,  enters the game. 
\claim{}  
To pursue our program we use the basis of $\bw\Vcal\otimes \bw\Vcal^*=\bigoplus_{k,l\geq 0}\bw^k\Vcal\otimes \bw^l\Vcal^*$ obtained as the  union of those induced on  $\bw^k\Vcal\otimes \bw^l\Vcal^*$ by $\bfb$ and $\bbeta$,  for all $k,l\geq 0$.  This is quite  straightforward, up to getting  aware of one main combinatorial point, i.e. that they  are best parametrised by the set 
$\ovPc$ of what, in Definition~\ref{def:bPartition}, lacking of a better terminology,  we called {\em bilateral partitions}.
More precisely, given $r\geq 0$,   we shall understand by  $\ovPc_r$ the set of all $r$-tuples $
\blamb=(\lambda_1,\ldots, \lambda_r)\subseteq\ZZ^r$, such that $\lambda_1\geq\cdots\geq 
\lambda_r$. We so have 
$$
\bw^k\Vcal=\bigoplus_{\bmu\in\ovPc_k}\QQ\wb^k_\bmu\qquad \mathrm{and}\qquad
\bw^l\Vcal^*=\bigoplus_{\bnu\in\ovPc_l}\QQ\wbet^l_\bmu,
$$
where 
$$
\wb^k_\bmu=b_{k-1+\mu_1}\w\cdots\w b_{\mu_k}\qquad \mathrm{and}\qquad 
\wbet^l_\bnu=\beta_{l-1+\nu_1}\w\cdots\w \beta_{\nu_l}.
$$
Then
$$
\Ecal(\bfz_k,\bfw_l^{-1})=\sum_{\bmu,\bnu\in\ovPc_k\otimes \ovPc_l}\wb^k_\bmu\otimes \wbet^l_\bnu\bfs_\bmu(\bfz_k)\bfs_\bnu(\bfw_l^{-1}),
$$
is 
the generating function of the distinguished basis  $\wb^k_\bmu\otimes \wbet^l_\bnu$ of $\bw^k\Vcal\otimes\bw^l\Vcal^*$,
where $\bfz_k$ and $\bfw_l^{-1}$ are, respectively,  $k$-tuples  $(z_1,\ldots,z_k)$ and  $l$-tuples  $(w_1^{-1},\ldots, w_l^{-1})$ of formal variables. Abusing  notation, we have chosen to denote by the same  symbols $\bfs_\bmu(\bfz_k)$ and $\bfs_\bnu(\bfw_l^{-1})$   natural extensions of the classical Schur polynomials occurring in the theory of symmetric functions as in, e.g., \cite[Section 3]{MacDonald} and/or \cite[Section 2.2.]{Fulyoung}. The difference with the classical ones is  that they are symmetric {\em rational functions} and do coincide with the usual Schur symmetric polynomials whenever  $\blamb\in\Pcal_r=\ovPc\cap\NN^r$.
We are now in position to anticipate the statement of our main result.

\medskip
\noindent
{\bf Theorem \ref{thm:thm34}.} {\em The {\em (DJKM bosonic)} action of $\Ecal(\bfz_k,\bfw_l^{-1})$ on $B(\xi)$ is given by
\be
\Ecal(\bfz_k,\bfw_l)=\exp\left(\sum_{n\geq 1}{1\over n}p_n(\bfz_k^{-1})p_n(\bfw_l)\right)\Gamma(\bfz_k,\bfw_l),\label{eq0:mnthm}
\ee
\vspace{-5pt}
where
\begin{enumerate}
 \item[i)] the map  $\Gamma(\bfz_k,\bfw_l):B(\xi)
\sra B(\xi)\llb\bfz_k^{\pm 1}, \bfw_l^{\pm 1}\rrb$ is the {\em vertex operator}
\be
R(\bfz_k,\bfw_l^{-1})\hskip-1pt\exp\left(\sum_{n\geq 1}x_n(p_n(\bfz_k)-
p_n(\bfw_l))\right)\hskip-2pt\exp\hskip-2pt\left(\sum_{n\geq 1}
{p_n(\bfz_k^{-1})-p_n(\bfw_l^{-1})\over n}{\d\over d x_n}\right);\label{eq0:gfint}
\ee
\item[ii)] the map $R(\bfz_k,\bfw_l^{-1}): B(\xi)\llb\bfz_k,\bfw_l^{-1}\rrb\sra B(\xi)\llb\bfz_k,\bfw_l^{-1}\rrb$ is the unique $B\llb\bfz_k,\bfw_l^{-1}\rrb$-linear  extension of 
$$\xi^m\mapsto \xi^{m+k-l}\displaystyle{\prod_{\begin{scriptsize}\begin{matrix}1\leq i\leq k\cr
1\leq j\leq l\end{matrix}\end{scriptsize}}{z_i^{m-l+1}\over w_j^{m-l+1}}};$$
\item[iii)] the  expression  $p_n(\bfz_k^{\pm })$ and $p_n(\bfw_l^{\pm 1})$ denote the Newton 
powers sums symmetric polynomials, in the variables $\bfz_k^{\pm 1}$ and $
\bfw_l^{\pm 1}$, i.e.
more explicitly
$$
p_n(\bfz_k^{\pm 1}):=z_1^{\pm n}+\cdots +z_k^{\pm n}\qquad\mathrm{and}\qquad p_n(\bfw_l^{\pm 1}):=w_1^{\pm n}+\cdots +w_l^{\pm n}.
$$
\end{enumerate}
}

\noindent
The meaning of formula \eqref{eq0:gfint} is that if $P(\bfx,\xi)\in B(\xi)$ is any polynomial, then the ``multiplication'' of $\wb^k_\bmu\otimes\wb^l_\bnu$ 
 is the coefficient of $\bfs_\bmu(\bfz_k)\bfs_\bnu(\bfw_l^{-1})$ in the expansion $\Ecal(\bfz_k,\bfw_l^{-1})P(\bfx,\xi)$. This may seem tricky. However multiplying the resulting expression by the product of the Vandermonde $\Delta_0(\bfx_k)\Delta_0(\bfw_l^{-1})$, it is sufficient to consider the coefficient of the less intimidating monomial $z_k^{k-1+\mu_1}\cdots z_1^{\mu_k}\cdot w_1^{-l+1-\nu_1}\cdots w_l^{-\nu_k}$ suffices. 
 
 To end up, reading formula \eqref{eq0:gfint} for  $k=l=1$, putting $z_1=z$ and $w_1=w$, one has $\bfs_{(i)}(z)=z^i$  and $\bfs_{(j)}(w^{-1})=w^{-j}$,  for all $i,j\in\ZZ$. By  the definition of the logarithm of an invertible formal power series:
$$
\exp\left(\sum_{n\geq 1}{1\over n}{w_n\over z^n}\right)={1\over 1-\displaystyle{w\over z}}
$$
and the fact that, in this case,  $R(z,w^{-1})\xi^m=\xi^m\displaystyle{z^m\over w^m}$, equality \eqref{eq0:mnthm} simplifies into
\be
\Ecal(z,w^{_1})_{|B\xi^m}
={\displaystyle{z^m\over w^m}{1\over  1-\displaystyle{w\over z}}}
\exp\left(\sum_{n\geq 1}x_n(z^n-w^n)\right)
\exp\left(-\sum_{n\geq 1}{z^{-
n}-w^{-n}\over n}{\d \over \d x_n}\,\right),\label{eq0:djkmmm}
\ee
which is precisely the original DJKM formula for the bosonic 
representation of $gl(\Vcal)$ (see e.g. \cite[Proposition 5.2]{KacRaRoz} or \cite[Section 1]{jimbomiwa}. This may look surprising indeed, because  comparing \eqref{eq0:gfint} with \eqref{eq0:djkmmm}, it is apparent that  \eqref{eq0:gfint} can be obtained from the DJKM expression simply by replacing the variables 
$z,w$ in \eqref{eq0:djkmmm} by the power sums of the $k$ and $l$-tuples of indeterminates needed to write the appropriate generating functions. As in our previous references 
\cite{BeCoGaVi,gln,SDIWP}, we have  borrowed methods from the theory of Hasse-Schmidt 
derivation on a exterior algebra, like in the book \cite{HSDGA}. The similarity of DJKM formula with our \eqref{eq0:mnthm}, however, 
makes us wonder whether there is any other argument to deduce our Theorem \ref{thm:thm34}  bypassing our methods.

\claim{\bf Organisation of the paper.} In the first section we recall some more or less known pre-requisites. We revise, in particular, the construction of the fermionic 
Fock space following \cite[Section 5]{SDIWP} as well as how to  extend the Schubert derivation on it.  A little background on 
Schur polynomials, mainly following \cite{MacDonald} but also 
\cite[Lecture 6]{KacRaRoz}, is  included as well. Section \ref{sec:sec2} is devoted to carefully  define the 
 generating function of the basis elements of $\bw^k\Vcal\otimes
\bw^l\Vcal^*$, that  is best suited to describe the fermionic and bosonic 
representation of $gl(\bw\Vcal)$. In this same section 
the natural notion of bilateral partition is also introduced. It is reasonable to suspect it somewhere hidden  in 
some less known literature.  
Section \ref{sec:sec3} eventually concerns the statement and proof of our main theorem which, as announced,  supplies the expression of both the fermionic 
and the bosonic expression of $gl(\bw\Vcal)$. The two cases 
are  treated in a  unified way, reflecting the fact inspiring the
references \cite{HSDGA, gln,SDIWP} that there is a very little,  if not any at all, substantial 
difference between the two spaces. Indeed,  as explained in 
\cite{BeCoGaVi},  the vertex operators occurring in the representation 
theory of the Heisenberg algebra,  come naturally to life, exactly the 
same, already at the level of multivariate Schubert derivations on 
exterior algebras. With no serious need, at least for the focused
purposes of our research, to cross the walls to enter in the realm of the infinite wedge powers, as however we did in the present contribution.

\section{Background and notation}

\claim{} \label{sec:sec11} We shall deal with a $\QQ$-vector space $\Vcal:=\bigoplus_{i\in
\ZZ}\QQ\cdot b_i$ and its restricted dual $\Vcal^*:=\bigoplus_{i\in\ZZ}\QQ
\cdot \beta_j$,  where $\beta_j \in \hom_\QQ(\Vcal, \QQ)$  is the unique 
linear form such that $\beta_j (b_i)=\delta_{ji}$. 
The   generating series of the basis elements of $\Vcal$ and $\Vcal^*$ 
 are, respectively:
\be 
\bfb(z)= \sum_{i \in \ZZ}b_i z^i \in V\llb z^{-1},z\rrb\qquad \mathrm{and}
\qquad \bbeta(w^{-1})= \sum_{j \in \ZZ}\beta_j w^{-j} \in V^*\llb w,w^{-1}\rrb.
\ee

\claim{\bf Hasse-Schmidt Derivations on $\bw\Vcal$.} A map $\Dcal(z):\bw\Vcal\sra \bw\Vcal\llb z\rrb$ is said to be {\em Hasse-Schmidt} (HS) derivation on $\bw\Vcal$ if $\Dcal(z)(\bfu\w\bfv)=\Dcal(z)\bfu\w\Dcal(z)\bfv$, for all $\bfu,\bfv\in\bw\Vcal$. Write $\Dcal(z)$ in the form $\sum_{j\geq 0}D_jz^j$, with $D_j\in\End_\QQ(\bw\Vcal)$. Then $\Dcal(z)$ is invertible in $\End_\QQ(\bw\Vcal)\llb z\rrb$ if and only if $D_0$ is invertible In this case $\Dcal(z)$ is invertible and its inverse $\overline{\Dcal}(z)$ is a HS--derivation as well.

\claim{\bf Schubert derivations.}
Consider the shifts endomorphisms $\sigma_{\pm 1}\in gl(\bw \Vcal)$ given by $\sigma_{\pm 1}b_j=b_{j\pm 1}$. 
By \cite[Proposition 4.1.13]{HSDGA}, there exist unique HS derivations on $\sigma_{\pm}(z):\bw\Vcal\sra \bw\Vcal\llb z^{\pm 1}\rrb$ such that
$$
\sigma_{\pm}(z)b_j=\sum_{i\geq 0}b_{j\pm i}z^{\pm i}.
$$
Let us denote by $\ovsig_{\pm}(z)$ their inverses in $\bw\Vcal\llb z^{\pm 1}\rrb$. Restricted to $\Vcal$ they work as follows
\be
\ovsig_+(z)b_j=b_j-b_{j+1}z\qquad \mathrm{and}\qquad \ovsig_-(z)b_j=b_j-b_{j-1}z^{-1}.
\ee
They are called {\em Schubert derivations} in the references \cite{HSDGA,gln,SDIWP}.

\claim{\bf Fermionic Fock space.}\label{sec:FFS} We quickly summarise the definition of 
the fermionic Fock space borrowed from \cite{SDIWP}. Let $[\Vcal]$ be a 
copy of $\Vcal$ (framed by square bracket to distinguish by the original $
\Vcal$ itself). It is the $\QQ$-vector space with basis $(\wb_m)_{m\in\ZZ}$ 
Identify $[\Vcal]$ with a sub-module of the tensor product $\bw\Vcal
\otimes_\QQ[\Vcal]$ via the map $\wb_m\mapsto 1\otimes \wb_m$. Let $W$ be 
the $\bw\Vcal$--submodule  of $\bw\Vcal\otimes_\QQ[\Vcal]$ generated by all the 
expressions $\{b_m\otimes \wb_{m-1}-\wb_m, b_m\otimes \wb_m\}_{m\in\ZZ}$. 
In formulas:
\[
W:=\bw\Vcal\otimes \big(b_m\otimes \wb_{m-1}-\wb_m\big)+ \bw\Vcal\otimes\big( b_m
\otimes \wb_m\big).\]
\bclm{\bf Definition.} {\em The fermionic Fock space is the $\bw\Vcal$-
module
\be
\Fcal:=\Fcal(\Vcal):={\bw\Vcal\otimes_\QQ[\Vcal]\over W}.
\ee
}
\eclm
Let $\bw\Vcal\otimes_\QQ[\Vcal]\sra \Fcal$ be the canonical projection. 
The class  of 
$u\otimes \wb_m$ in $\Fcal$ will be denoted  $u\w \wb_m$. Thus the 
equalities 
$b_m\w \wb_m=0$ and $b_m\w\wb_{m-1}=\wb_m$ hold in $\Fcal$. For all $m\in \ZZ$ 
and $\blamb\in\Pcal$ let, by 
definition
$$
\wb_{m+\blamb}:=\bfb^r_{m+\blamb}\w \wb_{m-r}=b_{m+\lambda_1}\w b_{m-1+
\lambda_2}\w\cdots\w b_{m-r+1+\lambda_r}\w\wb_{m-r}
$$
where $r$ is any positive integer such that $\ell(\blamb)\leq r$, which 
implicitly defines $\bfb^r_{m+\blamb}$ as an element of $\bw^r\Vcal_{\geq 
m-r+1}$, where by $\Vcal_{\geq j}$ we understand $\bigoplus_{i\geq j}\QQ
\cdot b_i$.
It turns out that $\Fcal$ 
is a graded $\bw\Vcal$-module:
$$
\Fcal:=\bigoplus_{m\in\ZZ}\Fcal_m,
$$
where
\be 
\Fcal_m:=\bigoplus_{\blamb\in\Pcal}\QQ\wb_{m+\blamb}=\bigoplus_{r\geq 0}
\bigoplus_{\blamb\in\Pcal_r}\QQ\bfb^r_{m+\blamb}\w \wb_{m-r}.\label{eq:Fm}
\ee
is the {\em fermionic Fock space} of charge $m$  \cite[p. 36]{KacRaRoz}.

\bclm{\bf Proposition.} {\em
\begin{enumerate}
\item[i)] The equality  $b_j\w\wb_m=0$ holds for all $j
\leq m$;
\item[ii)] The image of the map $\bw^r\Vcal\otimes \Fcal_m\sra \Fcal$ given by $(\bfu,\bfv)\mapsto u\w v$ is contained in $\Fcal_{m+r}$.
\end{enumerate}

}
\eclm
\proof They are \cite[Proposition 4.4 and 4.5]{SDIWP}.
\qed

\claim{\bf Extending Schubert derivations to $\Fcal$.} We now extend the Schubert derivations, in principle only defined on $\bw\Vcal$, on $\Fcal$ according to \cite{SDIWP} to which we refer to for more details. 
First we define their action on elements of the form $\wb_m$ by setting:
$$
\ovsig_-(z)\wb_m=\sigma_-(z)\wb_m:=\wb_m,\qquad \sigma_+(z)\wb_m:=\sigma_+(z)b_m\w \wb_{m-1}
$$
{and}
$$
\ovsig_+(z)\wb_m:=\sum_{j\geq 0}\wb_{m+(1^j)}z^j\qquad  \qquad 
$$
where $(1^j)$ denotes the partition with $j$ parts equal to $1$.
Finally, we set
\be
\sigma_{\pm}(z)\wb_{m+\blamb}=\sigma_{\pm}(z)\bfb^r_{m+\blamb}\w\sigma_\pm(z)\wb_{m-r}\qquad \mathrm{and}\qquad \ovsig_{\pm}(z)\wb_{m+\blamb}=\ovsig_{\pm}(z)\bfb^r_{m+\blamb}\w\wb_{m-r}.
\ee
\bclm{\bf Proposition.} {\em For all $m\in \ZZ$, Giambelli's formula for the Schubert derivation $\ovsig_+(z)$ holds:
\be
\wb_{m+\blamb}=\det(\sigma_{\lambda_j-j+i})\wb_m
\ee
}
\eclm
\proof See \cite[Proposition 5.13]{SDIWP}.\qed

We introduce now  an operator on $\Fcal$ which, in a sense, plays the role of the determinant of the shift endomorphism $\sigma_1$. We denote it by $\xi$. We shall understand it as the unique algebra endomorphism of $\bw \Vcal$ such that $\xi\cdot b_j=b_{j+1}$. Being an algebra homomorphism implies that
$$
\xi\bfb_{m+\blamb}=\bfb_{m+1+\blamb}
$$

It is  clearly invertible. Its inverse $\xi^{-1}$ is such that $\xi^{-1}b_j=b_{j-1}$. Secondly, we extend it to $\Fcal$ as follows:
\be
\xi \wb_{m+\blamb}=\xi(\bfb^r_{m+\blamb})\w \wb_{m+1+\blamb},\label{eq1:mstrd1}
\ee
where $r$ is any integer greater than the length of the partition $\blamb$. It is trivial to check that such a definition does not depend on the choice of $r> \ell(\blamb)$. So for instance
$$
\xi^{m'}\wb_{m+\blamb}=\wb_{m+m'+\blamb}.
$$

\claim{\bf Bosonic Fock space.} Let $B:=\QQ[\bfx]$, the polynomial ring in infinitely many indeterminates $\bfx:=(x_1,x_2,\ldots)$. 
As  a $\QQ$--vector space it possesses a basis of Schur polynomials parametrised by the set  $\Pcal$ of all partitions.  Moreover, $(S_1(\bfx), S_2(\bfx),\ldots)$ generate $B$ as a $\QQ$-algebra,  because $S_i(\bfx)$ is a polynomial of degree $i$, for all $i\geq 0$.   If $\blamb\in\Pcal$ one sets
\be
S_\blamb(\bfx)=\det(S_{\lambda_j-j+i}(\bfx))
\ee
where the sequence $(S_1(\bfx),S_2(\bfx), \ldots)$ is defined by
\be
\sum_{j\in \ZZ}S_j(\bfx)z^j=\exp(\sum_{i\geq 1}x_iz^i).\label{eq1:defsi}
\ee

Let $B(\xi):=B\otimes_\QQ\QQ[\xi^{-1},\xi]$ be the $\QQ[\xi]$-algebra  of $B$-valued Laurent polynomials in $\xi$. We shall refer to $B(\xi)$ as  the {\em bosonic Fock space.} It follows that
$$
B(\xi)=\bigoplus_{\begin{scriptsize}\begin{matrix}m\in\ZZ,\cr \blamb\in\Pcal\end{matrix}\end{scriptsize}}\QQ\cdot \xi^mS_\blamb(\bfx)
$$
\claim{} The space $\Fcal$ can be endowed with a structure of free $B(\xi)$-module generated by $\wb_0$  of rank one  generated by $\wb_0$ such that
$
\xi^mS_\blamb(\bfx)\wb_0=\wb_{\blamb},
$ by simply declaring
\begin{eqnarray}
\xi^mS_i(\bfx)\wb_{\blamb}&:=&\sigma_i\wb_{m+\blamb}.\label{eq1:mstrd2}
\end{eqnarray}
In fact
\begin{center}
\begin{tabular}{rcllr}
$\wb_{m+\blamb}$&$=$&$\xi^m\wb_\blamb$&& (Equation \eqref{eq1:mstrd1})\cr\cr
&$=$&$\xi^m\det(\sigma_{\lambda_j-j+i})\wb_0$&&\hskip 45pt (Giambelli's formula for Schubert derivations)\cr\cr
&$=$&$\xi^m\det(S_{\lambda_j-j+i})\wb_0$&&(by equality \eqref{eq1:mstrd2})\cr\cr
&$=$&$\xi^mS_\blamb(\bfx)\wb_0$&&(Definition of $S_\blamb(\bfx)$).
\end{tabular}
\end{center}
Equality \eqref{eq1:mstrd2}  can be also phrased by saying that $S_i(\bfx)$ is an eigenvalue of the $\QQ(\xi)$-linear map $\sigma_i:\Fcal\sra \Fcal$ with $\Fcal_m$ as eigenspaces. It 
 implies that
\be
\sigma_+(z)\wb_{m+\blamb}=\exp\left(\sum_{i\geq 1}x_iz^i\right)\wb_{m+\blamb},
\ee
i.e., abusing terminology,  $\exp(\sum_{i\geq 1}{x_iz^i})$ is an eigenvalue of $\sigma_+(z)$.

\bclm{\bf Lemma.}\label{lem:lem112} {\em

\begin{enumerate}
\item[i)] The Schubert derivations $\sigma_\pm(z),\ovsig_\pm(z)$ commute with multiplication by $\xi$, i.e.
\be
\xi\sigma_{\pm}(z)=\sigma_{\pm}(z)\xi\qquad \mathrm{and}\qquad \xi\ovsig_{\pm}(z)=\ovsig_{\pm}(z)\xi;\label{eq1:commr}
\ee
\item[ii)] by regarding the Schubert derivation $\sigma_-(z)$ (resp. $\ovsig_-(z)$)  as a map $B\sra B[z^{-1}]$ by setting $(\sigma_-(z)S_\blamb(\bfx))\wb_m=\sigma_-(z)\wb_{m+\blamb}$ (resp. $(\ovsig_-(z)S_\blamb(\bfx))\wb_m=\ovsig_-(z)\wb_{m+\blamb}$, one has:
\begin{eqnarray}
\ovsig_-(z)S_i(\bfx)&=&S_i(\bfx)-\displaystyle{S_{i-1}(\bfx)\over z}\label{eq1:23}\\ \cr
\sigma_-(z)S_i(\bfx)&=&\sum_{j=0}^i{S_{i-j}(\bfx)\over z^j};\label{eq1:24}
\end{eqnarray}
\item[iii)] the maps $\sigma_-(z)$ and $\ovsig_-(z)$ are $\QQ(\xi)$-algebra endomorphism of $B(\xi)$. In particular
\be
\sigma_-(z)S_\blamb(\bfx)=\det(\sigma_-(z)S_{\lambda_j-j+i}(\bfx))
\label{eq1:a27x}
\ee
and
\be 
 \ovsig_-(z)S_\blamb(\bfx)=\det(\ovsig_-(z)S_{\lambda_j-j+i}(\bfx));
 \label{eq1:a28x}
\ee
\item[iv)] the maps  $\sigma_-(z)$ and $\ovsig_-(z)$ act on $B$ as 
exponential of a first order differential operators, namely:
\be 
\sigma_-(z)S_\blamb(\bfx)=\exp\left(\displaystyle{\sum_{n\geq 1}{1\over 
nz^n}{\d\over \d x_n}}\right)S_\blamb(\bfx) \label{eq1:a28}
\ee
and
\be 
\ovsig_-(z)S_\blamb(\bfx)=\exp\left(-\displaystyle{\sum_{n\geq 1}{1\over 
nz^n}{\d\over \d x_n}}\right)S_\blamb(\bfx).\label{eq1:a27}
\ee

\end{enumerate}
}
\eclm
\proof
i)  First we show that the commutation holds on the exterior algebra $\bw 
\Vcal$. This is nearly obvious, because 
$$
\sigma_\pm(z)\xi b_{j}=\sigma_\pm(z)b_{j+1}=\sum_{i\geq 0}b_{j+1\pm i}
z^{\pm i}=\xi\sum_{i\geq 0}b_{j\pm i}z^{\pm i}=\xi\sigma_\pm(z)b_j
$$
The same holds for $\ovsig_{\pm}(z)$. We have
$$
\ovsig_\pm(z)\xi b_{j}=\ovsig_\pm(z)b_{j+1}=b_{j+1}-b_{j+1\pm 1} z^{\pm 1}
=\xi(b_{j}-b_{j\pm 1} z^{\pm 1})=\xi\,\ovsig_\pm(z)b_j.
$$ Secondly, the commutation rules hold for elements of the form $\wb_m$. 
In fact:
\begin{center}
\begin{tabular}{rcllr}
$\sigma_-(z)\xi\wb_m$&$=$&$\sigma_-(z)\wb_{m+1}$&&(Definition of $\xi$)\cr
\cr
&$=$&$\wb_{m+1}$&&($\sigma_-(z)$ acts as the identity)\cr\cr
&$=$&$\xi\wb_m=\xi\sigma_-(z)\wb_m$&&(Definition of $\xi$ and $\sigma_-(z)
$ acts\\
&&&& as the identity on $\wb_m$)
\end{tabular}
\end{center}
Similarly one sees that $\ovsig_-(z)\xi=\xi\ovsig_-(z)$.
The check for $\sigma_+(z)$ and $\ovsig_+(z)$ works analogously as 
follows. 
\begin{center}
\begin{tabular}{rcllr}
$\sigma_+(z)\xi\wb_m$&$=$&$\sigma_+(z)\wb_{m+1}$&&(Definition of $\xi$)\cr
\cr
&$=$&$\sigma_+(z)b_{m+1}\w \wb_{m}$&&(Definition of $\sigma_+(z)\wb_m$) 
\cr\cr
&$=$&$\sum_{i\geq 0} b_{m+1+i}z^i\w \wb_{m}$&&\hskip30pt (Definition of $\sigma_+
(z)b_m$)\cr\cr
&$=$&$\sum_{i\geq 0}\xi b_{m+i}\w \xi\wb_{m-1}=\xi\sigma_+(z)\wb_m$
\end{tabular}
\end{center}
and
\begin{center}
\begin{tabular}{rclr}
$\ovsig_+(z)\xi\wb_m$&$=$&$\ovsig_+(z)\wb_{m+1}$&(Definition of $\xi$)\cr
\cr
&$=$&$\sum_{j\geq 0}(-1)^{j}b_{m+1+(1^j)}\w \wb_{m-j}z^j$&(Definition of 
$\ovsig_+(z)\wb_{m+1}$) \cr\cr
&$=$&$\sum_{j\geq 0}(-1)^{j}\xi b_{m+(1^j)}\w \xi \wb_{m-1-j}z^j
$&\hskip-5pt(Definition of multiplying by $\xi$)\cr\cr
&$=$&$\xi\sum_{j\geq 0}(-1)^{j} b_{m+(1^j)}\w  \wb_{m-1-j}z^j=\xi\ovsig_+
(z)\wb_m$\cr\cr

\end{tabular}
\end{center}

\noindent
Let us show now that \eqref{eq1:commr} holds when evaluated 
against a general element of $\Fcal$. We  check for $\sigma_+
(z)$, the others being analogous and even easier. Let $\blamb$ be  any 
partition and $r$ any integer such that $\ell(\blamb)<r$. Then:

\begin{center}
\begin{tabular}{rcllr}
$\sigma_\pm(z)(\xi\wb_{m+\blamb})$&$=$&$\sigma_\pm(z)\wb_{m+1+\blamb}$&(definition of multiplication by $\xi$)\cr
\cr
&$=$&$\sigma_\pm(z) (\bfb^r_{m+1+\blamb}\w\wb_{m+1-r})$&(decomposition of $\wb_{m+1+\blamb}$)\cr\cr
&$=$&$\sigma_\pm(z)\bfb^r_{m+1+\blamb}\w\sigma_\pm(z)\wb_{m+1-r}$&($\sigma_{\pm}(z)$ is a derivation)\cr\cr
&$=$&$\sigma_\pm(z)\xi\bfb^r_{m+\blamb}\w\sigma_\pm(z)\xi\wb_{m-r}$&(definition of multiplication by $\xi$)\cr\cr
&$=$&$\xi\sigma_\pm(z)\bfb^r_{m+\blamb}\w\xi \sigma_\pm(z)\wb_{m-r}$&(Lemma \ref{lem:lem112}, item i))\cr\cr 
&$=$&$\xi\sigma_{\pm}(z)\wb_{m+\blamb}$.
\end{tabular}
\end{center}
The proof for the Schubert derivations $\sigma_-(z)$ and $\ovsig_\pm(z)$ works the same.

\medskip
\noindent
ii) The proof of this second statement works verbatim as in 
\cite[Proposition 5.3]{pluckercone}, where the  $S_i(\bfx)$ are denoted by 
$h_i$;

\medskip
\noindent
iii) In this case the check follows  by combining  \cite[Proposition 
7.1]{pluckercone} and \cite[Corollary 7.3]{pluckercone};

\medskip
\noindent
 iv) Recall  that $B(\xi)=\QQ(\xi)[S_1(\bfx),S_2(\bfx),\ldots]$.  Equation~\eqref{eq1:defsi} implies that
$$
{\d S_i(\bfx)\over \d x_j}=S_{i-j}(\bfx),
$$
Then \eqref{eq1:23}, e.g., says that
\be
\ovsig_-(z)S_i(\bfx)=\left(1-{1\over z}{\d\over \d x_1}\right)S_i(\bfx)=
\exp\left(-\sum_{n\geq 1}{1\over nz^n}{\d^n\over \d x_1^n}\right)S_i(\bfx)
\ee
Now $\displaystyle{\d^n\over \d x_1^n}S_i(\bfx)=\displaystyle{\d\over \d 
x_n}S_i(\bfx)$. Since $S_i(\bfx)$ generate $B$ as a $\QQ$-algebra 
and $\ovsig_-(z)$ are algebra homomorphisms coinciding on generators,  
\eqref{eq1:a27} follows. The proof of \eqref{eq1:a28} is analogous, but it also
follows from inverting both members of the equality \eqref{eq1:a27}, obtaining
$$
\sigma_-(z)=\exp\left(\sum_{n\geq 1}{1\over nz^n}{\d\over \d x_n}\right)
$$
\qed

\claim{} In the sequel we will need the following observation. Suppose that $\phi$ is any of the endomorphism $\sigma_{\pm i}$ of $\ovsig_{\pm j}$, for $i$ and $j$ arbitrary non negative integers and that
$$
\phi\wb_{m+\blamb}=\sum_{\bmu}a_\bmu\wb_{m+\bmu}.
$$
Then, for any $m'\in\ZZ$, 
$$
\sum_{\bmu}a_\bmu\wb_{m+m'+\bmu}=\phi\wb_{m+m'+\blamb}.
$$
The proof is based on the definition of multiplication by $\xi$.
\begin{eqnarray*}
\sum_{\bmu}a_\bmu\wb_{m+m'+\bmu}&=&\sum_{\bmu}a_\bmu\xi^{m'}\wb_{m+\bmu}
=\xi^{m'}\sum_{\bmu}a_\bmu\wb_{m+\bmu}\cr\cr
&=&\xi^{m'}\phi\wb_{m+\blamb}=\phi\xi^{m'}\wb_{m+\blamb}=\phi\wb_{m+m'+\blamb}.
\end{eqnarray*}

\section {The generating functions of the bases of $\bw^k\Vcal$ and $\bw^l
\Vcal^*$}\label{sec:sec2}

Let  $\bw \Vcal=\bigoplus_{k\geq 0}\bw^k\Vcal
$ and $\bw\Vcal^*=\bigoplus_{l\geq 0}\bw^l\Vcal^*$ be the exterior algebra of $\Vcal$ and $\Vcal^*$ respectively.
To describe the  bases of $\bw^k \Vcal$ and  $\bw^l \Vcal^*$ induced by the basis $\bfb$  of $\Vcal$ and of $\bbeta$ of $\Vcal^*$ (Cf. Section \ref{sec:sec11}),  we need to explain what we shall mean by {\em bilateral partition}. 
\bclm{\bf Definition.}\label{def:bPartition}
{\em 
A {\em bilateral partition} of length at most $r\geq 1$ is an element of the 
set:
$$
\ovPc_r:=\left\{\blamb:=(\lambda_1,\lambda_2,\dots,\lambda_r)\in \ZZ^r\,|
\, \lambda_1\geq\lambda_2\geq\dots\geq \lambda_r\right\}.
$$
}
\eclm
Clearly, $\Pcal_r:=\ovPc_r\cap \NN^r$ is the set of the usual partitions of length at most $r$, namely the non--increasing sequences of non--negative integers with at most $r$ non zero parts.
If $i_1>\cdots>i_k$ is a decreeasing sequence of integers, there exists 
one and only one bilateral partition $\bmu\in\ovPc_k$ such that
$i_j=k-j+\mu_j$. Therefore  $(\wb^k_\bmu)_{\bmu\in\ovPc_k}$ and $(\wbet^l_\bnu)_{\bnu\in
\ovPc_l}$ where:
 $$
 \wb^k_\bmu =b_{k-1+\mu_1} \w \dots \w b_{\mu_k}\qquad \mathrm{and}\qquad 
 \wbet^l_\bnu =\beta_{l-1+\nu_1} \w \dots \w \beta_{\nu_l},
 $$
 are $\QQ$-bases of $\bw^k\Vcal$ and $\bw^l\Vcal^*$ respectively.
Let $\bfz_k:=(z_1,\ldots, z_k)$ and $\bfw_k^{-1}:=(w_1^{-1},\ldots, 
w_k^{-1})$  be two ordered finite sequences of formal variables. 
The $\bw^k\Vcal$-valued formal power series
$$
\bfb(z_k)\w\cdots\w \bfb(z_1)
$$
vanishes whenever $z_i=z_j$, for all $1\leq i<j\leq k$. Therefore it is 
divisible by the Vandermonde determinant $\Delta_0(\bfz_k)=\prod_{1<\leq 
i<j\leq k}(z_j-z_i)$. We then define, for all $\blamb\in\ovPc_k$, the 
extended Schur polynomial
$$
\bfs_\blamb(\bfz_k)
$$
through the equality
\be
\sum_{\bmu\in \ovPc}\wb^k_\bmu\bfs_\bmu(\bfz_k)\Delta_0(\bfz_k):=\bfb(z_k)\w\cdots\w \bfb(z_1),\label{eq1:gfbz}
\ee
and therefore the expression
\be
\wb^k(\bfz_k):=\sum_{\bmu\in \ovPc_k}\wb^k_\bmu\bfs_\bmu(\bfz_k)\label{eq1:gfbet}
\ee
is a generating function of the basis elements of $\bw^k\Vcal$ induced by 
the given basis $\bfb$ of $\Vcal$. Similarly, a generating function for 
the basis elements $(\wbet^l_{\bnu})_{\bnu\in\ovPc_l}$ is given by
\be
\wbet^l(\bfw_l^{-1}):=\sum_{\bnu\in\ovPc_\bnu}\wbet^l_\bnu\cdot\bfs_
\bnu(\bfw_l^{-1}),\label{eq1:eqgfbett}
\ee
where $\bfs_\bnu(\bfw_l^{-1})$ is now defined, for all $\bnu\in\ovPc_l$, 
via the equality
\be
\sum_{\bnu\in \ovPc}\wbet^l_\bnu\bfs_\bnu(\bfw_l^{-1})
\Delta_0(\bfw_l^{-1}):=\bbeta(w_1^{-1})\w\cdots\w \bbeta(w_l^{-1}),
\label{eq1:eqgfbet}
\ee
 where
\be
 \Delta_0(\bfw^{-1}_l)=\prod_{1<\leq i<j\leq l}(w_j^{-1}-
 w_i^{-1})={\prod_{1\leq i<j\leq l}(w_i-w_j)\over \prod_{i=1}^lw_i^{l-1}}.
\ee
Notice the different numbering adopted for the variables $\bfz$ (formula \eqref{eq1:gfbz}) and  the 
variables $\bfw^{-1}$ (formula \eqref{eq1:eqgfbet})

\claim{\bf Remark.} 
If  $\blamb \subseteq \NN ^k$, then $\bfs_{\blamb}(\bfz_k)$ 
is 
the usual Schur symmetric polynomial in $(z_1, z_2 , \ldots ,z_k)$. If $
\blamb:=(\lambda_1,\lambda_2,\dots,\lambda_k)\in \ovPc_{k}$, with all $
\lambda_i < 0$, then
\be
\bfs_\blamb(\bfz_k)={\bfs_{-\blamb}(\bfz_k^{-1})\over z_1^{k-1}\cdots 
z_k^{k-1}}.
\ee
where $-\blamb=(-\lambda_k,-\lambda_{k-1},\ldots,-\lambda_1)$.
If $\lambda_1 > 0$ and $\lambda_k < 0$, instead
\be
\bfs_{\blamb_k}(\bfz_k)={\bfs_{(\lambda_1 + \lambda_k, \cdots ,
\lambda_{k-1} 
+ \lambda_k , 0)}(\bfz_k)\over \prod_{j=0}^{k}z_j^{\lambda_k}}.
\ee
It is then clear that all $\bfs_{\blamb}(\bfz)$, where $\blamb$ runs on $
\Pcal_k$, are $\QQ$-linearly independent. The same holds true for $\Delta_0(\bfw_l^{-1})$.
\claim{} Let $\beta\in\Vcal^*$. The contraction $\beta\lrcorner:\bw\Vcal\sra \bw\Vcal$ can be depicted via the following diagram:
\be
\begin{vmatrix}\beta(b_{r-1+\lambda_1})&\beta(b_{r-2+\lambda_2}) &\ldots&
\beta(b_{\lambda_r})\cr\cr
b_{r-1+\lambda_1}&b_{r-2+\lambda_2}&\ldots&b_{\lambda_r}
\end{vmatrix}\label{eq2:contb}
\ee
to be read as follows. The scalar $\beta(b_{r-j+\lambda_j})$ is the 
coefficient of the element of $\bw^{r-1}\Vcal$ obtained by removing the $j
$-th exterior factor from $\wb^r_\blamb$. 

The contraction of $\bw^r\Vcal$ 
against $\wbet^l_\bnu\in \bw^l\Vcal^*$ is well defined as well.  It is an element of $\bw^{r-l}
\Vcal$ which can be  represented as (See 
\cite{BeCoGaVi}):
\be
\wbet^l_\bnu\lrcorner \wb^r_\blamb=\begin{vmatrix}\beta_{l-1+\nu_1}
(b_{r-1+\lambda_1})&\ldots&\beta_{l-1+\nu_1}(b_{\lambda_r})\cr\vdots&
\ddots&\vdots\cr
\beta_{\nu_l}(b_{r-1+\lambda_1})&\ldots&\beta_{\nu_l}(b_{\lambda_r})\cr\cr
b_{r-1+\lambda_1}&\ldots&b_{\lambda_r}\label{eq2:contra}
\end{vmatrix}
\ee
to be read as follows. The Laplace-like expansion of the array 
\eqref{eq2:contra} along the first row is an alternating linear 
combination of contractions of elements of $\bw^{k-1}\Vcal$ against 
elements of $\bw^{l-1}\Vcal^*$. Having already set the case $k=1$ in \eqref{eq2:contb}, we have described it completely.

\claim{} Although it may be easily guessed, let us now  make precise the 
definition of the contraction of an element of $\Fcal$ against an element 
of $\bw^l\Vcal^*$. Giving  the definition on bases elements $
\wb_{m+\blamb}$ of $\Fcal$ and  $\wbet^l_\bnu$  ($\bnu:=(\nu_1\geq\ldots
\geq\nu_l)$ of $\bw^l\Vcal^*$ will suffice.  Let $r\geq 0$ such that $\ell(\blamb)\leq 
r$ and $\nu_l\geq m-r$ and define:
$$
\wbet^l_{\bnu}\lrcorner \wb_{m+\blamb}:=(\wbet^l_\bnu\lrcorner \wb^r_{m+\blamb})\w\wb_{m-r}.
$$ 
It is straightforward to see that the definition does not depend on the choice of the non-negative integer $r>\ell(\blamb)$.

\claim{} Let 
\be
\Ecal(\bfz_k,\bfw_l^{-1})=\wb^k(\bfz_k)\otimes \wbet^l(\bfw_l^{-1})=
\sum_{\bmu,\bnu}\wb^k_\bmu\otimes \wbet^l_\bnu \bfs_\bmu(\bfz_k)\bfs_
\bnu(\bfw_l^{-1}),\label{eq2:genfbk}
\ee
be the generating function of the basis of $\bw^k\Vcal\otimes \bw^l\Vcal^*$. It defines two maps  
\begin{eqnarray}
\Ecal_f(\bfz_k,\bfw_l^{-1}):\Fcal\sra \Fcal\llb\bfz_k,\bfw_l, \bfz_k^{-1},\bfw_l^{-1}]
\end{eqnarray}
and
\begin{eqnarray}
\Ecal_b(\bfz_k,\bfw_l^{-1}):=B(\xi)\sra B(\xi)\llb\bfz_k,\bfw_l, \bfz_k^{-1},\bfw_l^{-1}]
\end{eqnarray}
which we distinguish by putting a subscript in the notation and satisfying the compatibility relation imposed by the boson-fermion correspondence.  More precisely we define:
\be
\Ecal_f(\bfz_k,\bfw_l^{-1})\wb_{m+\blamb}
:=\wb^k(\bfz_k)\w \wbet^l(\bfw_l^{-1})\lrcorner \wb_{m+\blamb}\label{eq2:repdf}
\ee
and 
\be
\big(\Ecal_b(\bfz_k,\bfw_l^{-1})\xi^mS_\blamb(\bfx)\big)\wb_0=\Ecal_f(\bfz_k,\bfw_l^{-1})\wb_{m+\blamb}\label{eq2:repdb}
\ee
where we have used  the notation of \eqref{eq1:gfbet} and \eqref{eq1:eqgfbet}.
\claim{\bf Products of Schubert derivations.}  To further elaborate the shape of \eqref{eq2:repdf} and \eqref{eq2:repdb}, we need to introduce 
the following new piece of notation.
Let
\be
\sigma_+(\bfz_k)=\sigma_+(z_1)\cdots\sigma_+(z_k),\qquad \qquad\ovsig_+
(\bfz_k)=\ovsig_+(z_1)\cdots\ovsig_+(z_k),\label{eq2:prodz}
\ee
and
\be
\sigma_-(\bfw_l)=\sigma_-(w_1)\cdots\sigma_-(w_l),\qquad \qquad\ovsig_
-(\bfz_l)=\ovsig_-(w_1)\cdots\ovsig_-(w_l).\label{eq2:prodw}
\ee
Equalities \eqref{eq2:prodz} and \eqref{eq2:prodw} must be read in $\End_
\QQ(\bw\Vcal)\llb \bfz_k\rrb$ and  $\End_\QQ(\bw\Vcal)\llb \bfw_l^{-1}\rrb
$ respectively. They are {\em multivariate HS-derivations} of $\bw \Vcal$ 
in the following sense: i) they are {\em multi-variate} because  are $
\End_\QQ(\bw\Vcal)$ 
formal power series in more than one indeterminate, namely  $\bfz_k:=(z_1,
\ldots,z_k)$ and $\bfw_l^{-1}:=(w_1^{-1},\ldots,w_l^{-1})$, and ii) are {\em HS derivations}, being 
compatible with the wedge product:
$$
\sigma_\pm(\bfz_k)(\bfu\w \bfv)=\sigma_\pm(\bfz_k)\bfu\w \sigma_
\pm(\bfz_k)\bfv\qquad \mathrm{and}\qquad \ovsig_\pm(\bfz_k)(\bfu\w \bfv)=
\ovsig_\pm(\bfz_k)\bfu\w \ovsig_\pm(\bfz_k)\bfv.
$$

\bclm{\bf Proposition.}\label{prop:propwbet} {\em The following equality holds:
\be
\bbeta(w_1^{-1})\w\cdots\w \bbeta(w_l^{-1})\lrcorner \wb_{m+\blamb}={\Delta_0(\bfw_l^{-1})\over \prod_{j=1}^lw_j^{m-l+1}}\ovsig_+(\bfw_l)\sigma_-(\bfw_l)\wb_{m-l+\blamb}.
\ee
}
\eclm 
\proof
If $l=1$ formula \eqref{eq2:bwb} reads as
$$
\beta(w_1^{-1})\lrcorner \wb_{m+\blamb}=w_1^{-m}\ovsig_+(w_1)\sigma_-(w_1)\wb_{m-1+\blamb}
$$  
and this is precisely \cite[Proposition 6.13]{SDIWP}. Assume the formula holds for $l-1\geq 0$. For notational simplicity let  $\bfw_l\ssetminus w_1:=(w_2,\ldots,w_l)$ and $\bfw_l^{-1}\ssetminus w_1^{-1}:=(w_2^{-1},\ldots, w_l^{-1})$. Then
\begin{center}
	\begin{tabular}{rcllll}
		&&$\bbeta(w_1^{-1})\w\cdots\w \bbeta(w_l^{-1})\lrcorner \wb_{m+\blamb}$&\cr\cr
		&$=$&$\bbeta(w_1^{-1})\lrcorner \left(\bbeta(w_2^{-1})\cdots\w \bbeta(w_l^{-1})\lrcorner \wb_{m+\blamb}\right)$
		\hskip45pt (Associativity of "$\w $" )\cr\cr
		&$=$&$\bbeta(w_1^{-1})\lrcorner
		\left(\displaystyle{\Delta_0 (\bfw_l^{-1}\ssetminus w_1^{-1})\over\prod_{j=2}^lw_j^{m-l+2}}\right)\ovsig_+(\bfw_l\ssetminus w_1)\sigma_-(\bfw_l\ssetminus w_1)\wb_{m-l+1+\blamb}
		$\cr\cr
		&$=$&$w_1^{-m+l-1}\displaystyle{\Delta_0 (\bfw_l^{-1}\ssetminus w_1^{-1})\over  \prod_{j=2}^lw_j^{m-l+2}}
		\ovsig_+(w_1)\sigma_-(w_1)
		\ovsig_+(\bfw_k\ssetminus w_1)\sigma_-(\bfw_k\ssetminus w_1)
		\wb_{m-k+\blamb}$\cr
	\end{tabular}
\end{center}		
Now we use the commutation rule:
\begin{eqnarray*}
	\sigma_-(w_1)\ovsig_+(\bfw_k\ssetminus w_1)&=&
	\ovsig_+(\bfw_k\ssetminus w_1)\sigma_-(w_1)
	\left(1-{w_2 \over w_1}
	\right)\cdots 
	\left(1-{w_l \over w_1}
	\right)
\end{eqnarray*}
From which
\begin{center}
	\begin{tabular}{rclrr}		
		&$=$&$\displaystyle{w_1^{-m+l-1} \over w_1^{l-1}}(w_1-w_2)\cdots 
		(w_1-w_l)
		\displaystyle{\Delta_0 (\bfw_l^{-1}\ssetminus w_1^{-1})\over \prod 
		w_j^{m+l-2}}\cdot\ovsig_+(\bfw_l)\sigma_-(\bfw_l)
		\wb_{m-k+\blamb}$\cr\cr
		&$=$&$\displaystyle{w_1^{-m}\over w_1^{l-1}w_2\cdots w_l}\prod
		\left(\displaystyle{1\over w_j}-\displaystyle{1\over w_1}\right)
		\displaystyle{\Delta_0 (\bfw_l^{-1}\ssetminus w_1^{-1})\over \prod 
		w_j^{m+l-2}}\ovsig_+(\bfw_l)\sigma_-(\bfw_l)
		\wb_{m-k+\blamb}$
		\cr\cr

		&$=$&$\displaystyle{\Delta_0(\bfw_l^{-1})\over \prod_{j=1}
		^lw_j^{m-l+1}}\ovsig_+(\bfw_l)\sigma_-(\bfw_l)\wb_{m-k+\blamb}
		$\cr
	\end{tabular}
\end{center}
as desired. \qed

\bclm{\bf Corollary.}\label{prop:prop24} {\em The generating function \eqref{eq1:eqgfbett} acts on $\Fcal$ according to:
\be
\sum_{\bnu\in\ovPc_l}\wbet^l_\bnu \bfs_\bnu(\bfw_l^{-1})\lrcorner\wb_{m+
\blamb}=\prod_{j=1}^lw_j^{-m+l-1}\ovsig_+(\bfw_l)\sigma_-(\bfw_l)\wb_{m-l+
\blamb}.\label{eq2:bwb}
\ee
}

\proof
It is a consequence of equality \eqref{eq1:eqgfbet} and of Proposition \ref{prop:propwbet} up to  dividing by the Vandermonde determinant. \qed

\eclm
\bclm{\bf Proposition.} \label{prop:propgfb}{\em For all $k\geq 1$:
$$
\bfb(z_k)\w\cdots\w \bfb(z_1)\w \wb_{m+\blamb}=\prod_{j=1}^k z_j^{m+1}\Delta_0(\bfz_k)\sigma_+(\bfz_k)\ovsig_-(\bfz_k^{-1})\wb_{m+k+\blamb}.
$$
}
\eclm
\proof By induction on $k \geq 1.$ If $k=1$, the formula reads as 
$$
\bfb(z_1)\w \wb_{m+\blamb}=z_1^{m+1}
\sigma_+(\bfz_1)\ovsig_-(\bfz_1)\wb_{m+1+\blamb}
$$
and this is Proposition $6.9$ in \cite{SDIWP}. Assume the formula holds for $ k-1 \geq 0 $. Then,
\begin{center}
	\begin{tabular}{rclrr}
		&&$\bfb(z_k)\w\cdots\w \bfb(z_1)\w \wb_{m+\blamb}$&\cr\cr
		&$=$&$\bfb(z_k)\w 
		\left(\bfb(z_{k-1}) \w\cdots\w \bfb(z_1)\w \wb_{m+\blamb}\right)$
		\hskip40pt (Associativity of "$\w $" )\cr\cr
		&$=$&$\bfb(z_k) \w z_{k-1}^{m+1} \cdots z_1^{m+1} 
		\sigma_+(\bfz_{k-1})\ovsig_-(\bfz_{k-1})\wb_{m+k-1+\blamb}\cdot \Delta_0(\bfz_{k-1})$\cr\cr
		&$=$&$z_k^{m+k}z_{k-1}^{m+1} \cdots z_1^{m+1}\Delta_0(\bfz_{k-1})\sigma_+(\bfz_k)\ovsig_-(\bfz_k)\sigma_+(\bfz_{k-1})\ovsig_-(\bfz_{k-1})\wb_{m+k+\blamb}$\cr\cr
		&$=$&$z_k^{m+k+1}\prod_{j=1}^{k-1}z_j^{m+1} \prod_{j=1}^{k-1}
		\displaystyle{\left(1-{z_j \over z_k}\right)}  \Delta_0(\bfz_{k-1})\cdot$\cr\cr
		&&$\cdot\,\,\sigma_+(z_1)\sigma_+(\bfz_{k-1})
		\ovsig_-(z_1)
		\ovsig_-(\bfz_{k-1})\wb_{m+k-1+\blamb}$\cr\cr
		&$=$&$\displaystyle{z_k^{m+k+1}\over z_k^{k-1}}\prod_{j=1}^{k-1}z_j^{m+1} \prod_{j=1}^{k-1}
		\displaystyle{\left(z_k-z_j\right)}  \Delta_0(\bfz_{k-1})\cdot\sigma_+(\bfz_{k})\ovsig_-(\bfz_{k})\wb_{m+k-1+\blamb}$\cr\cr
		&$=$&$\displaystyle{\prod_{j=1}^k}z_j^{m+1}
		\Delta_0(\bfz_k)
		\sigma_+(\bfz_k)
		\ovsig_-(\bfz_k)\wb_{m+k+\blamb}$\cr
	\end{tabular}
\end{center}
as desired.
\qed
\bclm{\bf Corollary.}\label{cor:genbk} {\em The generating function \eqref{eq1:gfbz} acts on on the basis element $\wb_{m+\blamb}\in \Fcal$according to:
\be
\sum_{\bmu\in\ovPc_k}\wb^k_\bmu\bfs_\bmu(\bfz_k)\w \wb_{m+\blamb}=\displaystyle{\prod_{j=1}^k}z_j^{m+1}
		\sigma_+(\bfz_k)
		\ovsig_-(\bfz_k)\wb_{m+k+\blamb}
\ee
}
\eclm
\proof
By Proposition \ref{prop:propgfb}, using expression \eqref{eq1:gfbz}, dividing by the Vandermonde $\Delta_{0}(\bfz_k)$.\qed

\section{Fermionic and Bosonic Vertex Representation of $gl(\bw \Vcal)$.}\label{sec:sec3}

\bclm{\bf Lemma.}\label{prop:comgs} {\em The following commutation rules holds in $\End_\QQ(\Fcal)[z^{-1},w\rrb$
\begin{eqnarray}
\ovsig_-(z)\ovsig_+(w)&=&\left(1-{w\over z}\right)^{-1}\ovsig_+(w)\ovsig_-(z)\label{eq:comm2}\\ \cr
&=&\exp\left(\sum_{n\geq 0}{1\over n}{w^n\over z^n}\right)\ovsig_+(w)\ovsig_-(z)\label{eq:comm2b}
\end{eqnarray}
}
\eclm
\proof Formula~\eqref{eq:comm2} is  \cite[Proposition 8.4, Formula (54)]{SDIWP} and \eqref{eq:comm2b} uses the equality of formal power series $(1-x)^{-1}=\exp(\sum_{n\geq 1}x^n/n)$. \qed

\bclm{\bf Proposition.}\label{prop:prop32} {\em Let $\displaystyle{p_n(\bfz_k^{-1})=\sum_{i=1}^kz_i^{-n}}$ and $\displaystyle{p_n(\bfw_l)=\sum_{j=1}^lw_j^{n}}$ (the symmetric power sums Newton polynomials). The following equalities holds on $\End_{\QQ(\xi)}B(\xi)$:
\be
\ovsig_-(\bfz_k)=\prod_{j=1}^k\ovsig_-(z_j)=\exp\left(-\sum_{n\geq 1}{1\over n}p_n(\bfz_k^{-1}){\d\over \d x_n}\right).\label{eq3:for41}
\ee
and
\be
\sigma_-(\bfw_l)=\prod_{j=1}^l\sigma_-(z_j)=\exp\left(\sum_{n\geq 1}{1\over n}p_n(\bfw_l^{-1}){\d\over \d x_n}\right)\label{eq3:for42}
\ee
Therefore
\be
\ovsig_-(\bfz_k)\sigma_-(\bfw_l)\wb_{m+\blamb}=\left[\exp\left(-\sum_{n\geq 1}{1\over n}
(p_n(\bfz_k^{-1})-p_n(\bfw_l^{-1})){\d\over \d x_n}\right)\xi^mS_\blamb(\bfx)\right]\wb_0.\label{eq3:for43}
\ee
}
\eclm
\proof
The operators
$$
\sum_{n\geq 1}{1\over n}{1\over z_i^n}{\d\over \d x_n},\quad  \sum_{n\geq 
1}{1\over n}{1\over z_j^n}{\d\over \d x_n},\quad \sum_{n\geq 1}{1\over n}
{1\over w_p^n}{\d\over \d x_n}, \quad \sum_{n\geq 1}{1\over n}{1\over 
w_q^n}{\d\over \d x_n}
$$
commute for all choices of $1\leq i,j\leq k$ and $1\leq p,q\leq l$. Then 
the product of their exponential is the exponentials of their sum:
\begin{eqnarray*}
\ovsig_-(\bfz_k)=\prod_{j=1}^k\ovsig_-(z_j)&=&\prod_{j=1}^k\exp\left(-
\sum_{n\geq }{1\over nz^n_j}{\d\over \d x_n}\right)\cr\cr
&=&\exp\left(-\sum_{n\geq 1}{1\over n}\left({1\over z_1^n}+\cdots+{1\over 
z_k^n}\right){\d\over \d x_n}\right)\cr\cr
&=&\exp\left(-\sum_{n\geq 1}{1\over n}p_n(\bfz_k^{-1}){\d\over \d x_n}
\right),
\end{eqnarray*}
which validates \eqref{eq3:for41}. Formula \eqref{eq3:for42} is 
checked analogously.  Formula~\eqref{eq3:for43} follows from 
\eqref{eq3:for41} and \eqref{eq3:for42} and using again the fact that the 
operators $\displaystyle{\sum_{n\geq 1}{1\over n}p_n(\bfz_k^{-1}){\d\over \d x_n}}$ and 
$\displaystyle{\sum_{n\geq 1}{1\over n}p_n(\bfw_l^{-1}){\d\over \d x_n}}$ commute. Thus:
\begin{eqnarray}
\ovsig_-(\bfz_k)\sigma_-(\bfw_l)&=&\exp\left(-\sum_{n\geq 1}{1\over n}
p_n(\bfz_k^{-1}){\d\over \d x_n}\right)\exp\left(\sum_{n\geq 1}{1\over n}
p_n(\bfw_l^{-1}){\d\over \d x_n}\right)\cr\cr
&=&\exp\left(-\sum_{n\geq 1}{1\over n}(
p_n(\bfz_k^{-1})-p_n(\bfw_l^{-1})){\d\over \d x_n}\right).
\end{eqnarray}
\qed
\bclm{\bf Proposition.}\label{prop:commsk} {\em The following commutation rules holds:
\be
\ovsig_-(\bfz_k)\ovsig_+(\bfw_l)=\exp\left(\sum_{n\geq 1}{1\over n}p_n(\bfw_l)p_n(\bfz_k^{-1})\right)\ovsig_+(\bfw_l)\ovsig_-(\bfz_k)\label{eq3:forcor32}.
\ee
}
\eclm
\proof
We first prove that
\be
\ovsig_-(\bfz_k)\ovsig_+(\bfw_l)=\prod_{i=1}^k\prod_{j=1}^l\left(1-{w_j\over z_i}\right)^{-1}\ovsig_+(\bfw_l)\ovsig_-(\bfz_k)\label{eq3:prefor}
\ee
For $k=l=1$ the formula is Proposition~\eqref{prop:comgs}. Suppose it holds for $k-1\geq 1$ and $l=1$. Then
\begin{center}
\begin{tabular}{rcll}
$\ovsig_-(\bfz_k)\ovsig_+(\bfw_1)$&$=$&$\displaystyle{\prod_{i=1}^k\ovsig_-(z_i)}\cdot\ovsig_+(w_1)$&\hskip-10pt (definition of $\ovsig_+(\bfz_k)$)\cr\cr
&$=$&$\hskip-10pt \left(1-\displaystyle{w_1\over z_k}\right)^{-1}\displaystyle{\prod_{i=1}^{k-1}}\ovsig_-(z_i)\ovsig_+(w_1)\ovsig_-(z_k)$&\hskip-34pt(first step of induction on $l$)\cr\cr
&$=$&$\hskip-10pt \left(1-\displaystyle{w_1\over z_k}\right)^{-1}\displaystyle{\prod_{i=1}^{k-1}}\left(1-\displaystyle{w_1\over z_i}\right)^{-1}\ovsig_+(w_1)\displaystyle{\prod_{i=1}^{k-1}}\ovsig_-(z_i) \ovsig_-(z_k)$ &\hskip-3pt (inductive hypothesis\cr &&& \hskip 7pt on $k$)\cr\cr
&$=$&$\displaystyle{\prod_{i=1}^{k}}\left(1-\displaystyle{w_1\over z_i}\right)^{-1}\ovsig_+(w_1)\ovsig_-(\bfz_k)$ &\hskip-3pt (definition of $\ovsig_-(\bfz_k)$).
\end{tabular}
\end{center}
Suppose now that \eqref{eq3:prefor} holds for all $k\geq 1$ and $l-1\geq 0$. Then
\begin{center}
\begin{tabular}{rcll}
$\ovsig_-(\bfz_k)\ovsig_+(\bfw_l)$&$=$&$\ovsig_-(\bfz_k)\cdot\ovsig_+(w_l)
\ovsig_-(\bfw_{l-1})$&\cr\cr
&$=$&$\displaystyle{\prod_{i=1}^k}\left(1-\displaystyle{w_l\over z_i}
\right)^{-1}\ovsig_+(w_l)\ovsig_-(\bfz_k)\ovsig_+(\bfw_{l-1})$&\cr\cr
&$=$&$\displaystyle{\prod_{j=1}^l}\left(1-\displaystyle{w_l\over z_i}
\right)^{-1} \prod_{\begin{scriptsize}\begin{matrix}{{1\leq i\leq k}}\cr{ 
1\leq j\leq l-1}\end{matrix}\end{scriptsize}}\left(1-\displaystyle{w_j
\over z_i}\right)^{-1}\ovsig_+(w_l)\ovsig_+(\bfw_{l-1})\ovsig_-(\bfz_k)$&
\cr\cr
&$=$&$\displaystyle{\prod_{\begin{scriptsize}\begin{matrix}{{1\leq i\leq 
k}}\cr{ 1\leq j\leq l}\end{matrix}\end{scriptsize}}}\left(1-
\displaystyle{w_j\over z_i}\right)^{-1}\ovsig_+(\bfw_l)\ovsig_-(\bfz_k)$,&
\end{tabular}
\end{center}
which is precisely \eqref{eq3:prefor}.  To  phrase \eqref{eq3:prefor} in 
the form \eqref{eq3:forcor32} one first notice that
$$
\left(1-{w_j\over z_i}\right)^{-1}=\exp\left(\sum_{n\geq 1}{1\over n}
{w_j^n\over z_i^n}\right).
$$
By a simple manipulation one sees that
$$
\displaystyle{\prod_{\begin{scriptsize}\begin{matrix}{{1\leq i\leq k}}\cr{ 
1\leq j\leq l}\end{matrix}\end{scriptsize}}}\left(1-\displaystyle{w_j\over 
z_i}\right)^{-1}=\prod_{\begin{scriptsize}\begin{matrix}{{1\leq i\leq k}}
\cr{ 1\leq j\leq l}\end{matrix}\end{scriptsize}}\exp\left(\sum_{n\geq 1}
{1\over n}{w_j^n\over z_i^n}\right)=\exp\left(\sum_{n\geq 1}{1\over n}
p_n(\bfw_l)p_n(\bfz_k^{-1})\right)
$$
as desired. \qed

\claim{} Let $R_f(\bfz_k,\bfw_l^{-1}):\Fcal\sra \Fcal[\bfz_k^{\pm 1},\bfw_l^{\pm 1}]$ defined on homogeneous elements as:
$$
R_f(\bfz_k,\bfw_l^{-1})\wb_{m+\blamb}={\prod_{i=1}^kz_i^{m-l+1}\over 
\prod_{j=1}^lw_j^{m-l+1}}\xi^{k-l}\wb_{m+\blamb}
$$
and $R_b(\bfz_k,\bfw_l^{-1})\in \hom_{\QQ[\xi]}(B(\xi),[\bfz_k^{\pm 1},\bfw_l^{\pm 1}])$ defined by
$$
(R_b(\bfz_k,\bfw_l^{-1})\xi^mS_\blamb(\bfx))\wb_0=R_f(\bfz_k,\bfw_l^{-1})\wb_{m+\blamb}
$$
from which
$$
R_b(\bfz_k,\bfw_l^{-1})\cdot 1={\prod_{i=1}^kz_i^{m-l+1}\over 
\prod_{j=1}^lw_j^{m-l+1}}\xi^{k-l}
$$
\bclm{\bf Proposition.} {\em
The map $R_f(\bfz_k,\bfw_l^{-1})$ commutes with Schubert 
derivations, in the sense that
$$
\sigma_{\pm}(\bfz_k)R_f(\bfz_k,\bfw_l^{-1})=R_f(\bfz_k,\bfw_l^{-1})
\sigma_{\pm}(\bfz_k)\qquad \mathrm{and}\qquad \ovsig_{\pm}(\bfz_k)R_f(\bfz_k,
\bfw_l^{-1})=R_f(\bfz_k,\bfw_l^{-1})\ovsig_{\pm}(\bfz_k).
$$
}
\eclm
\proof
It is enough to prove that it commutes with $\sigma_{\pm i}$ and $\ovsig_{\pm j}$, $i,j\geq 0$, which are by definition $\QQ[\bfx_k,\bfw_l^{-1}]$-linear. First of all recall that the product  $\sigma_{\pm i}\wb_{m+\blamb}$ ($\blamb\in\Pcal_r$) is ruled by some Pieri's formulas
$$
\sigma_{\pm i}\wb_{m+\blamb}=\sum_{\bmu\in P_{\pm}}\wb_{m+\bmu},
$$
where $P_+$ (resp. $P_-$) is the set of all partitions $\mu_1\geq 
\mu_2\geq\cdots\geq \mu_r$ ($r\geq \ell(\blamb)$) such that
$\mu_1\geq \lambda_1\geq\cdots\geq\mu_k\geq\lambda_k$ and $|\bmu|=|\blamb|
+i$ (resp. $\lambda_1\geq \mu_1\geq\lambda_2\geq \mu_2\geq\cdots\geq
\lambda_r\geq  \mu_r$ and $|\bmu|=|\blamb|-i$). 
Then we have
\begin{eqnarray*}
\sigma_{\pm i}R_f(\bfz_k,\bfw_l^{-1})\wb_{m+\blamb}&=&\sigma_{\pm i}
\prod_{\begin{scriptsize}\begin{matrix}1\leq i\leq k\cr1\leq j\leq l 
\end{matrix}\end{scriptsize}}
{z_i^{m+l-1}\over w_j^{m-l+1}}\,\,\xi^{k-l}\wb_{m+\blamb}
=\prod_{\begin{scriptsize}\begin{matrix}1\leq i\leq k\cr1\leq j\leq l 
\end{matrix}\end{scriptsize}}{z_i^{m+l-1}\over w_j^{m-l+1}}\xi^{k-l}
\sigma_{\pm i}\wb_{m+
\blamb}\cr\cr\cr
&=&
\prod_{\begin{scriptsize}\begin{matrix}1\leq i\leq k\cr1\leq j\leq l 
\end{matrix}\end{scriptsize}}{z_i^{m+l-1}\over w_j^{m-l+1}}\xi^{k-l}
\sum_{\bmu\in P_{\pm}}\wb_{m+\bmu}=R_f(\bfz_k,\bfw_l^{-1})\sum_{\bmu\in P_
\pm}\wb_{m+\bmu}\cr\cr
&=&R_f(\bfz_k,\bfw_l^{-1})\sigma_{\pm i}\wb_{m+\blamb}
\end{eqnarray*}
Thus $\sigma_{\pm }(\bfz_k)$ commutes with $R_f(\bfz_k,\bfz_l^{-1})$ and 
so do $\ovsig_{\pm}(\bfz_k)$. Indeed:
\begin{eqnarray*}
\ovsig_{\pm}(\bfz_k)R_f(\bfz_k,\bfz_l^{-1})&=&\ovsig_{\pm}(\bfz_k)R(\bfz_k,
\bfz_l^{-1})\sigma_{\pm}(\bfz_k)\ovsig_{\pm}(\bfz_k)\cr\cr
&=&\ovsig_\pm(\bfz_k)\sigma_\pm(\bfz_k)R(\bfz_k,\bfw_l^{-1})\ovsig_
\pm(\bfz_k)\cr\cr
&=&R_f(\bfz_k,\bfw_l^{-1})\ovsig_\pm(\bfz_k).\hskip240pt \qed
\end{eqnarray*}

\bclm{\bf Theorem.}\label{thm:thm34} {\em Notation as in \eqref{eq2:repdf} 
and \eqref{eq2:repdb}.
Then:
\be
\Ecal_f(\bfz_k,\bfw_l^{-1})=\exp\left(\sum_{n\geq 1}
{1\over n}p_n(\bfw_l)p_n(\bfz_k^{-1})\right)\Gamma_f(\bfz_k,\bfw_l)
\label{eq3:fermnth}
\ee
and
\be
\Ecal_b(\bfz_k,\bfw_l^{-1})=\exp\left(\sum_{n\geq 
1}
{1\over n}p_n(\bfw_l)p_n(\bfz_k^{-1})\right)\Gamma_b(\bfz_k,\bfw_l.)\label{eq3:ecalb}
\ee
where the fermionic and bosonic vertex operators are, respectively
\begin{eqnarray}
\Gamma_f(\bfz_k,\bfw_l)&=&R_f(\bfz_k,\bfw_l^{-1})\sigma_+(\bfz_k)\ovsig_+
(\bfw_l)\ovsig_-(\bfz_k)\sigma_-(\bfw_l)\cr\cr
&=&R_f(\bfz_k,\bfw_l^{-1})\exp\left(\sum_{n\geq 1}x_n(p_n(\bfz_k)-p_n(\bfw_l))\right)\ovsig_-(\bfz_k)\sigma_-(\bfw_l).\label{eq3:fvoiif}
\end{eqnarray}
and
\be
\hskip-6pt \Gamma_b(\bfz_k,\bfw_l)\hskip-2pt=\hskip-2pt R_b(\bfz_k,\bfw_l^{-1})\exp\hskip-2pt\left(\sum_{n\geq 1}x_n(p_n(\bfz_k)-p_n(\bfw_l))\right)\hskip-2pt\exp\hskip-2pt\left(\hskip-3pt-\sum_{n\geq 1}
{p_n(\bfz_k^{-1})-p_n(\bfw_l^{-1})\over n}{\d\over \d x_n}\right)\label{eq3:bvof}
\ee
}
\eclm
\proof
We have:
\begin{center}
\begin{tabular}{rcll}
$\Ecal_f(\bfz_k,\bfw_l^{-1}))\wb_{m+\blamb}$&$=$&$\wb^k(\bfz_k)\w 
\wbet^l(\bfw_l^{-1})\lrcorner \wb_{m+\blamb}$&\hskip-48pt (definition of $
\Ecal_f(\bfz_k,\bfw_l^{-1}))$\cr\cr
&$=$&$\wb^k(\bfz_k)\w \displaystyle{\prod_{j=1}^l}w_j^{-m+l-1}\ovsig_+(\bfw_l)\sigma_-
(\bfw_l^{-1})\wb_{m-l+\blamb}$&(Corollary~\ref{prop:prop24})\cr\cr
&$=$&$\displaystyle{\prod_{i=1}^kz_i^{m-l+1}\over \prod_{j=1}^lw_j^{m-l
+1}}\sigma_+(\bfz_k)\ovsig_-(\bfz_k)\ovsig_+(\bfw_l)\sigma_-(\bfw_l^{-1})
\wb_{m+k-l+\blamb}$&(Corollary~\ref{cor:genbk})\cr\cr
&$=$&$R(\bfz_k,\bfw_l^{-1})\sigma_+(\bfz_k)\ovsig_-(\bfz_k)\ovsig_+
(\bfw_l)\sigma_-(\bfw_l^{-1})\wb_{m+\blamb}$&(Definition \cr
&&&of $R(\bfz_k,\bfw_l^{-1})$)
\end{tabular}
\end{center}
By invoking the commutation relation proven in Proposition 
\ref{prop:commsk}, one obtains
\be
\Ecal_f(\bfz_k,\bfw_l^{-1})\wb_{m+\blamb}=\exp
\left(\displaystyle{\sum_{n\geq 1}}\displaystyle{1\over n}
p_n(\bfw_l)p_n(\bfz_k^{-1})\right)R(\bfz_k,\bfw_l^{-\1})\sigma_+(\bfz_k)
\ovsig_+(\bfw_l)\ovsig_-
(\bfz_k)\sigma_-(\bfw_l)\wb_{m+\blamb}\label{eq3:ferpic}
\ee
which already prove that the expression of $\Ecal_f(\bfz_k,\bfw_l^{-1})$ 
is precisely \eqref{eq3:fermnth}.
To continue with, the $B(\xi)$-module structure of $\Fcal$ says that $
\Fcal_m$ is an 
eigenspace of $\sigma_+(\bfz_k)\ovsig_+(\bfw_l)$ with eigenvalue
\begin{eqnarray}
\prod_{i=1}^k\exp\left(\sum_{n\geq 1}x_nz_i^n\right)\prod_{j=1}^l\exp
\left(-\sum_{n\geq 1}x_nw_j^n\right)&=&\exp\left(\sum_{n\geq 1}
x_np_n(\bfz_k)\right)\exp\left(-\sum_{n\geq 1}x_np_n(\bfw_l)\right)\cr\cr
\cr
&=&\exp\left(\sum_{n\geq 1}x_n(p_n(\bfz_k)-p_n(\bfw_l))\right).
\label{eq3:eignv}
\end{eqnarray}
Thus formula \eqref{eq3:ferpic}, up to replacing $\sigma_+(\bfz_k)\ovsig_
+
(\bfw_l)$ by its eigenvalue \eqref{eq3:eignv} with respect to $\Fcal$, is 
precisely \eqref{eq3:fermnth} with $\Gamma_f(\bfz_k,\bfw_l)$ given by 
expression \eqref{eq3:fvoiif}.
To prove \eqref{eq3:ecalb} we recall that
\begin{eqnarray*}
(\Ecal_b(\bfz_k,\bfw_l^{-1})\xi^mS_\blamb(\bfx))\wb_0&=&\Ecal_f(\bfz_k,
\bfw_l^{-1})\wb_{m+\blamb}=\exp\left(\displaystyle{\sum_{n\geq 1}}
\displaystyle{1\over n}
p_n(\bfw_l)p_n(\bfz_k^{-1})\right)\Gamma_f(\bfz_k,\bfw_l)\wb_{m+\blamb}
\end{eqnarray*}
Now
\begin{eqnarray*}
\Gamma_f(\bfz_k,\bfw_l)\wb_{m+\blamb}&=&R_f(\bfz_k,\bfw_l^{-1})\exp\left(\sum_{n
\geq 1}x_n(p_n(\bfz_k)-p_n(\bfw_l))\right)\ovsig_-(\bfz_k)\sigma_-(\bfw_l)
\wb_{m+\blamb}
\end{eqnarray*}
However, by Proposition \ref{prop:prop32}, formula~\eqref{eq3:for43}, 

$$
\ovsig_-(\bfz_k)\sigma_-(\bfw_l)\wb_{m+\blamb}=\left[\exp\left(-\sum_{n
\geq 1}
{p_n(\bfz_k^{-1})-p_n(\bfw_l^{-1})\over n}{\d\over \d x_n}\right)\xi^mS_
\blamb(\bfx)\right]\wb_0.
$$
which shows that 
$$\Gamma_f(\bfz_k,\bfw_l)\wb_{m+\blamb}=(\Gamma_b(\bfz_k,\bfw_l^{-1})
\xi^mS_\blamb(\bfx))\wb_0
$$
proving the theorem.\qed

\claim{\bf Remark.} In formula \eqref{eq3:ecalb} let us set $k=l=1$ and 
call $z=z_1$ and $w=w_1$. Then $p_n(z^{\pm 1})=z^{\pm n}$ and $p_n(w^{\pm 
1})=w^{\pm n}$.
Then
$$
\Ecal_b(z,w)=\exp\left(\sum_{n\geq 1}{1\over n}{w^n\over z^n}\right)
\Gamma_b(z,w)
$$
where
$$
\Gamma_b(z,w)=R_b(z,w)\exp(\sum_{n\geq 1}x_n(z^n-w^n))\exp\left(-\sum_{n
\geq 1}{z^{-n}-w^{-n}\over n}{\d\over \d x_n}\right)
$$
Keeping into account that 
$$
\exp\left(\sum_{n\geq 1}{1\over n}{w^n\over z^n}\right)={1\over 1-\displaystyle{w\over 
z}}
$$ and using the definition of $R_b(z,w^{-1})$ one sees that
$$
\Ecal_b(z,w)_{|B^{(m)}}={\displaystyle{z^m\over w^m}\over 1-\displaystyle{w\over z}}\exp(\sum_{n\geq 1}x_n(z^n-w^n))\exp\left(-\sum_{n\geq 1}{z^{-n}-w^{-n}\over n}{\d\over \d x_n}\right)
$$
which is the celebrated DJKM formula.

\smallskip
\noindent{\bf Acknowledgments.} This work is one of the topics touched in the Ph.D. thesis of the first author, mostly redacted during her hosting at the Department of Mathematical Sciences of Politecnico of Torino under the sponsorship of Ministry of Science of the Islamic Republic of Iran. The second author profited of the support of Finanziamento Diffuso della Ricerca (no. 53$\_$RBA17GATLET)  and Progetto di Eccellenza del Dipartimento di Scienze Matematiche, 2018--2022, no.
E11G18000350001. The project also benefitted the partial 
support of INDAM-GNSAGA, PRIN ``Geometria delle Variet\`a Algebriche''.

For discussions and criticisms we want to 
 primarily thank Inna Scherbak, who, as in \cite{BeCoGaVi}, first suggested us to generalise the DJKM picture in the way as now stands in the present paper and Parham Salehyan for many useful redactional hints. Finally, we are indebted to Joachim Kock and Andrea T. Ricolfi for  their carefully reading and many other kinds of assistance.

\bibliographystyle{amsplain}

\medskip
\medskip

\parbox[t]{3in}{{\rm Ommolbanin~Behzad}\\ 
	\smallskip \vspace{-15pt}
	
	{\tt \href{mailto:behzad@iasbs.ac.ir}{behzad@iasbs.ac.ir}}\\
	{\it Department of Mathematics\\
		Institute for Advanced Studies\\ in Basic Sciences (IASBS)\\ P.O.Box 45195-1159  Zanjan\\ IRAN}}
	\parbox[t]{3in}{{\rm Letterio Gatto}\\
	{\tt \href{mailto:letterio.gatto@polito.it
		}{letterio.gatto@polito.it}}\\
	{\it Dipartimento di Scienze Mate--}\\
	{matiche, Politecnico di Torino}\\
	{\it ITALY}}

\end{document}